\documentclass[11pt,twoside,a4paper]{article}
\usepackage{amssymb}
\usepackage{amsmath}
\usepackage[dvips]{graphicx}
\usepackage[all,2cell,ps]{xy}
\title{A Generalized Macaulay Theorem and Generalized Face Rings}
\author{Eran Nevo \footnote{Institute of Mathematics, Hebrew
University, Jerusalem, Israel, E-mail address: eranevo@math.huji.ac.il. }}
 \newtheorem{thm}{Theorem}[section]
 
  \newtheorem{cor}[thm]{Corollary}
 \newtheorem{lem}[thm]{Lemma}
\newtheorem{de}[thm]{Definition}
 \newtheorem{prop}[thm]{Proposition}
\newtheorem{prob}[thm]{Problem}

\begin{document}
\maketitle
\begin{abstract}
We prove that the $f$-vector of members in a certain class of meet
semi-lattices satisfies Macaulay inequalities $0\leq
\partial^{k}(f_{k})\leq f_{k-1}$ for all $k\geq 0$.
We construct a large family of meet semi-lattices belonging to
this class, which includes all posets of multicomplexes, as well
as meet semi-lattices with the "diamond property", discussed by
Wegner \cite{W}, as special cases. Specializing the proof to the
later family, one obtains the Kruskal-Katona inequalities and
their proof as in \cite{W}.

For geometric meet semi-lattices we construct an analogue of the
exterior face ring, generalizing the classic construction for
simplicial complexes. For a more general class, which also
includes multicomplexes, we construct an analogue of the
Stanley-Reisner ring. These two constructions provide algebraic
counterparts (and thus also algebraic proofs) of Kruskal-Katona's
and Macaulay's inequalities for these classes, respectively.
\end{abstract}

\section{Introduction}\label{Introduction}
Let us review the characterization of $f$-vectors of finite
simplicial complexes, known as the
Sch\"{u}tzenberger-Kruskal-Katona theorem (see \cite{Bol} for a
proof and for references). Let $C$ be a (finite, abstract)
simplicial complex, $f_{i}=|\{S\in C: |S|=i+1\}|$.
$f=(f_{-1},f_{0},...)$ is called the $f$-vector of $C$ (note that
$f_{-1}=1$). For any two integers $k,n>0$ there exists a unique
expansion
\begin{equation} \label{nk}
n= {n_{k}\choose k} + {n_{k-1}\choose k-1}+...+ {n_{i}\choose i}
\end{equation}
such that $n_{k}>n_{k-1}>...>n_{i}\geq i \geq 1$ (details in
\cite{Bol}). Define the function $\partial_{k-1}$ by
 $$\partial_{k-1}(n)=
{n_{k}\choose k-1} + {n_{k-1}\choose k-2}+...+ {n_{i}\choose
i-1},\ \ \partial _{k-1}(0)=0.$$
\begin{thm}[Sch\"{u}tzenberger-Kruskal-Katona] \label{KKthm}
$f$ is the $f$-vector of some
simplicial complex iff $f$ ultimately vanishes and
\begin{equation} \label{partial_{k-1}(n)}
\forall k\geq 0\ \   0\leq \partial_{k}(f_{k})\leq f_{k-1}.
\end{equation}
\end{thm}
For a ranked meet semi-lattice $P$, finite at every rank, let
$f_{i}$ be the number of elements with rank $i+1$ in $P$, and set
${\rm rank}(\hat{0})=0$ where $\hat{0}$ is the minimum of $P$. the
$f$-vector of $P$ is $(f_{-1},f_0,f_1,...)$.

$P$ has the $diamond\ property$ if for every $x,y\in P$ such that
$x<y$ and ${\rm rank} (y)-{\rm rank} (x)=2$ there exist at least
two elements in the open interval $(x,y)$. The closed interval is
denoted by $[x,y]=\{z\in P: x\leq z\leq y \}$.

We identify a simplicial complex with the poset of its faces
ordered by inclusion. The following generalization of Theorem
\ref{KKthm} is due to Wegner \cite{W}.

\begin{thm}[Wegner]\label{Wthm}
Let $P$ be a finite ranked meet semi-lattice with the diamond
property. Then its $f$-vector ultimately vanishes and satisfies
(\ref{partial_{k-1}(n)}).
\end{thm}
For $\hat{x}\in P$ define $P(\hat{x})=\{x\in P: \hat{x}\leq x\}$
and let $y'\prec  y$ denote $y$ covers $y'$.

\begin{lem}\label{cond* lemma}
For a ranked meet semi-lattice $P$, the diamond property is
equivalent to satisfying the following condition:

(\textbf{*}) For every $\hat{x}\in P$, $x$ which covers $\hat{x}$
and $y$ such that $y \in P(\hat{x})$ and $y\not =\hat{x}$, there
exists $y'\in P(\hat{x})$ such that $y'\prec y$ and $x\nleq y'$.
\end{lem}

A multicomplex (on a finite ground set) can be considered as an
order ideal of monomials $I$ (i.e. if $m|n \in I$ then also $m\in
I$) on a finite set of variables. Its $f$-vector is defined by
$f_{i}=|\{m\in I: {\rm deg} (m)=i+1\}|$ (again $f_{-1}=1$). Define
the function  $\partial^{k-1}$ by $$\partial^{k-1}(n)=
{n_{k}-1\choose k-1} + {n_{k-1}-1\choose k-2}+...+ {n_{i}-1\choose
i-1},\ \
\partial^{k-1}(0)=0,$$ w.r.t the expansion (\ref{nk}).
\begin{thm}[Macaulay,\cite{M}] \label{Mthm} (More proofs in \cite{CL,St})\label{M thm}
 $f$ is
the $f$-vector of some multicomplex iff $f_{-1}=1$ and
\begin{equation} \label{partial^{k-1}(n)}
\forall k\geq 0\ \    0\leq \partial^{k}(f_{k})\leq f_{k-1}.
\end{equation}
\end{thm}

\begin{figure}\label{Fig1}
\newcommand{\edge}[1]{\ar@{-}[#1]}
\newcommand{\lulab}[1]{\ar@{}[l]^<<{#1}}
\newcommand{\rulab}[1]{\ar@{}[r]^<<{#1}}
\newcommand{\ldlab}[1]{\ar@{}[l]^<<{#1}}
\newcommand{\rdlab}[1]{\ar@{}[r]_<<{#1}}
\newcommand{\node}{*+[O][F-]{ }}
\centerline{ \xymatrix{
 & & & & \circ \rulab{y} \edge{d} \edge{ddlll} & & & &\circ \rulab{y} \edge{d} \edge{dlll} &\\
& \circ \lulab{x_{i+1}} \edge{d}\ & & & \bullet \rulab{\exists y'} \edge{ddlll} & \circ \lulab{x_{r}} \edge{d} & & & \bullet \rulab{\exists y'} \edge{dlll} &\\
& \circ \lulab{x_{i}} \edge{d} & & & & \circ \lulab{x_{r-1}} \edge{dd}\\
& \circ \lulab{x_{i-1}} \edge{d}\\
& \circ \lulab{\hat{x}} & & & & \circ \lulab{\hat{x}} } } \caption
{The parallelogram property for $i<r$ (left) and for $i=r$
(right).}
\end{figure}
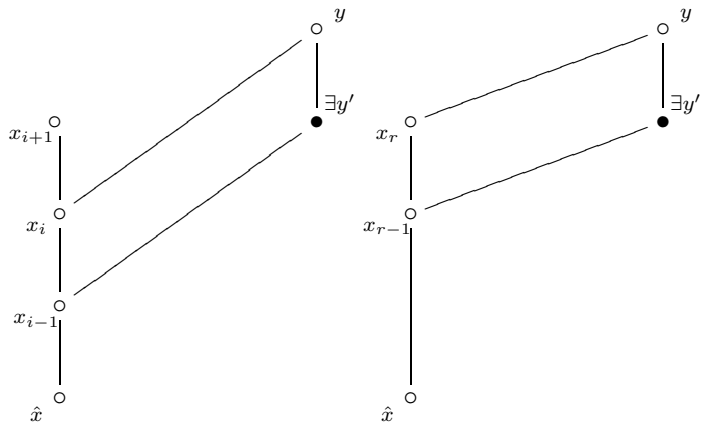

\begin{de}(Parallelogram property)\label{def**}
A ranked poset $P$ is said to have the $parallelogram$ $property$
if the following condition holds:

(\textbf{*}\textbf{*}) For every $\hat{x}\in P$ and $y\in
P(\hat{x})$ such that $y\not =\hat{x}$, if the chain
$\{\hat{x}=x_{0}\prec x_{1}\prec ...\prec x_{r} \}$ equals the
closed interval $[\hat{x},x_{r}]$ ($r>0$) and is maximal w.r.t.
inclusion such that $r<{\rm rank}(y)$ (the rank of $y$ in the
poset $P(\hat{x})$), and if $x_{i}<y$ and $x_{i+1}\nleqslant y$ for
some $0< i\leq r$, then there exists $y'\in P(\hat{x})$ such that
$y'\prec y$, $x_{i-1}<y'$ and $x_{i}\nleq y'$. For $i=r$ interpret
$x_{r+1}\nleqslant y$ as: $[\hat{x},y]$ is not a chain.
\end{de}
See Figure \ref{Fig1} for an illustration of the parallelogram
property. Note that condition (\textbf{*}) of Lemma \ref{cond*
lemma} implies condition (\textbf{*}\textbf{*}) of Definition
\ref{def**} (with $1$ being the only possible value of $r$), and
that posets of multicomplexes satisfy the parallelogram property.
\\
We identify a multicomplex with the poset of its monomials ordered
by division. We now generalize Theorem \ref{Mthm}.
\begin{thm} \label{THM}
Let $P$ be a ranked meet semi-lattice, finite at every rank, with
the parallelogram property. Then its $f$-vector satisfies
(\ref{partial^{k-1}(n)}) and $f_{-1}(P)=1$.
\end{thm}
For generalizations of Macaulay's theorem in a different direction
('compression'), see e.g. \cite{CL,WW}.
\\

In Section \ref{Proof} we prove Theorem \ref{THM} and construct a large family
 of meet semi-lattices satisfying its hypothesis.

Theorems \ref{KKthm} and \ref{Mthm} have algebraic counterparts in
terms of face rings and algebraic shifting. No such interpretation
is known for Theorems \ref{Wthm} and \ref{THM}. In Section
\ref{AlgShift} we extend Theorems \ref{KKthm} and \ref{Mthm} by
constructing analogues of the exterior and symmetric face rings,
respectively. More specifically, we define an exterior algebraic
shifting operation for geometric meet semi-lattices, and a
symmetric algebraic shifting operation for a common generalization
of geometric meet semi-lattices and multicomplexes.

\section{Macaulay inequalities}\label{Proof}
We provide proofs of Theorem \ref{THM} and Lemma \ref{cond*
lemma}, and construct a large class of examples for which Theorem
\ref{THM} applies.

$Proof\  of\  Theorem\  \ref{THM}$: Clearly $f_{-1}(P)=1$. Let us
show that $P$ satisfies the inequalities (\ref{partial^{k-1}(n)}).
Let $X$ be the set of rank $k+1$ ($\leq {\rm rank} (P)$) elements
in $P$, and denote its shadow by $\partial X$, i.e. $\partial
X=\{p\in P: \exists x\in X, p\prec  x \}$. We will show that
$|\partial X| \geq
\partial^{k}(|X|)$, which clearly proves Theorem \ref{THM}.

The proof is by induction on $k$ and on $|f_{k}|$. The case $k=0$
is trivial, as well as the case $|f_{k}|=1$ for any $k$. So assume
$k>0$.

Let us introduce some notation: Let $\hat{0}\not = x_{r}\in P$ be
such that the interval $[\hat{0},x_{r}]$ is maximal w.r.t.
inclusion such that it is a chain $\{\hat{0}=x_{0}\prec x_{1}\prec
...\prec  x_{r}\}$ and $x_{r}<x$ for some $x\in X$ (hence $r\leq
k$). For $0\leq i\leq r$, denote $P_{i}=\{p\in P: x_{i+1}\nleq p,
x_{i}< p \}$ and $X_{i}=X\cap P_{i}$. Thus $P=\uplus_{0\leq i\leq
r}P_{i}\uplus [\hat{0},x_{r}]$. In addition, $\partial
X=\uplus_{0\leq i\leq r}(\partial X\cap P_{i})$, unless $r=k$, in
which case $\{x_{k}\}$ should be added to that union. Let
$\hat{X_{i}}$ denote the elements of $X_{i}$ considered as
elements of the induced meet semi-lattice $P(x_{i})$. Thus,
$\partial (\hat{X_{i}})\subseteq
\partial X \cap P_{i}$ unless $x_{k}\in \partial (\hat{X_{i}})$, a
case in which $i=k$ and $\partial (\hat{X_{k}})=\{x_{k}\}$. Hence
\begin{equation} \label{partial hat}
   |\partial \hat{X_{i}}|\leq |\partial X \cap P_{i}| \ \ 0\leq i
   \leq \min \{r,k-1\},
\end{equation}
and for $r=k$ $|\partial \hat{X_{k}}|=1$. By the parallelogram
property, for any $0\leq i\leq \min \{r,k-1\}$ and $y\in X_{i+1}$,
there exists $y'\in \partial \{y\} \cap P_{i}$ (for $i=r \
X_{r+1}=\emptyset$). Note that $y'$ 's arising from different $y$
's are distinct: suppose $y'\in
\partial X \cap P_{i}$ arises from two different $y\in X_{i+1}$,
then as $P$ is a meet semi-lattice $x_{i+1}\leq y'$, a
contradiction. We deduce that
\begin{equation} \label{i+1}
   |X_{i+1}|\leq |\partial X \cap P_{i}| \ \ 0\leq i
   \leq \min \{r,k-1\}.
\end{equation}
Combining (\ref{partial hat}) and (\ref{i+1}) we get that
\begin{equation} \label{last}
    |\partial X | =
    \begin{cases}
        1+\sum_{0\leq i\leq k-1}|\partial X \cap P_{i}| \geq 1+\sum_{0\leq i\leq k-1}\max  \{|X_{i+1}|,|\partial \hat{X_{i}}|\}
        & \textrm{if $r=k$} \\ \\
        \sum_{0\leq i\leq r}|\partial X \cap P_{i}|  \geq \sum_{0\leq i\leq r}\max  \{|X_{i+1}|,|\partial \hat{X_{i}}|\} & \textrm{if $r<k$.}
    \end{cases}
\end{equation}
By induction hypothesis, $|\partial \hat{X_{i}}|\geq
\partial^{k-i}(|\hat{X_{i}}|) = \partial^{k-i}(|X_{i}|)$ for $0\leq i\leq \min \{r,k-1\}$ (the
induction on $k$ implies it for $i\neq 0$, and the induction on
$|f_{k}|$ implies it for $i=0$). We need the following simple
Lemma due to Bj\"{o}rner and Vre\'{c}ica: (One uses Theorem \ref{M
thm} to prove it.)
\begin{lem}(Lemma 3.2 of \cite{BV})\label{Lemma}
For $k>0$, the function $\partial^{k}$ satisfies for all
non-negative integers $n_{i}$ and $r<k$:
$$\partial^{k}(\sum_{0\leq i\leq r}n_{i})\leq \sum_{0\leq i\leq r}\max  \{n_{i+1},\partial^{k-i}(n_{i})\},$$
$$\partial^{k}(1+\sum_{0 \leq i \leq k}n_{i})\leq 1+\sum_{0\leq i\leq k-1}\max  \{n_{i+1},\partial^{k-i}(n_{i})\}.$$
\end{lem}
By Lemma \ref{Lemma} we see that both right hand sides of (\ref{last}) are $\geq
\partial^{k}(\sum_{0\leq j\leq r}|X_{j}|+\delta_{r,k})$, where $\delta_{r,k}$ is Kronecker's delta.
Using the fact that $\partial^{k}$ is nondecreasing, the right hand side of
(\ref{last}) is $\geq \partial^{k}(|X|)$ (as $X=\uplus_{0\leq i\leq r}X_{i}$ for $r<k$, and $|X|\leq 1+\sum_{0\leq i\leq r}|X_{i}|$ for $r=k$). Hence $|\partial X|\geq
\partial^{k}(|X|)$ as desired. $\square$
\\
\textbf{Remark}: If $P$ satisfies the diamond property, then
$X=X_{0}\uplus X_{1}$ and $|\partial X|\geq \sum_{0\leq i\leq
1}\max \{|X_{i+1}|,|\partial \hat{X_{i}}|\}$ (here
$X_{2}=\emptyset$), an inequality which implies the Kruskal-Katona
inequalities for $P$, via an inequality for the function
$\partial_{k}$, analogous to the one in Lemma \ref{Lemma}, which
is given in \cite{EW}. This is how the proof given in \cite{W}
argues.
\\
$Proof\  of\  Lemma\  \ref{cond* lemma}$: Condition (*) clearly
implies the diamond property. Conversely, we argue by induction on
$r=\rm{rank}(y)-\rm{rank}(\hat{x})$. For $r=1$, take $y'=\hat{x}$.
For $r=2$, this is the diamond property. For $r>2$, assume $x<y$
(otherwise the assertion is trivial). There exists $z$ such that
$x<z\prec y$. By the induction hypothesis, there exists $z'$ such
that $z'\prec z$ and $x\nleq z'$. By the diamond property applied
to the pair $(z',y)$, there exists $y'$ such that $z'\prec y'\prec
y$ and $y'\neq z$. Now, $x\nleq y'$ as otherwise we get $x\leq
y'\wedge z =z'$, a contradiction. $\square$

\textbf{Example F}: Let $(L,<,r)$ be a finite ranked meet
semi-lattice with partial order relation $<$ and rank function
$r$. Denote its minimum by $\hat{0}$. Associate with each
$\hat{0}\neq l\in L$ a collection $F(l)$ of multichains in the
interval $(\hat{0},l]$ which is closed w.r.t. the following
partial order on multichains in $L\setminus \{\hat{0}\}$: Let
$\underline{a}=(a_m\leq...\leq a_1\leq a_0)$ and
$\underline{b}=(b_k\leq...\leq b_1\leq b_0)$ be multichains in
$L\setminus \{\hat{0}\}$ and define $\underline{a}<'
\underline{b}$ iff $m\leq k$, $a_{i}\leq b_{i}$ for all $0\leq
i\leq m$ and $\underline{a}\neq \underline{b}$. $F(l)$ is closed
if $\underline{a}<'\underline{b}\in F(l)$ implies
$\underline{a}\in F(l)$.

\begin{figure}\label{FigL'}
\newcommand{\edge}[1]{\ar@{-}[#1]}
\newcommand{\lulab}[1]{\ar@{}[l]^<<{#1}}
\newcommand{\rulab}[1]{\ar@{}[r]^<<{#1}}
\newcommand{\rrulab}[1]{\ar@{}[rr]^<<{#1}}
\newcommand{\ldlab}[1]{\ar@{}[l]^<<{#1}}
\newcommand{\rdlab}[1]{\ar@{}[r]_<<{#1}}
\newcommand{\rrdlab}[1]{\ar@{}[rr]_<<{#1}}
\newcommand{\node}{*+[O][F-]{ }}
\centerline{ \xymatrix@C=2pt{
1,1,1\edge{dr}&1,14\edge{dr}\edge{d}&4,14\edge{d}\edge{drrrrrrrrr}&1,12\edge{dll}\edge{d}&2,12\edge{dl}\edge{dr}&2,2,2\edge{d}&
1234&2,23\edge{dl}\edge{dll}&3,23\edge{d}\edge{dll}&3,3,3\edge{dl}&3,34\edge{dl}\edge{dll}&4,34\edge{dll}\edge{d}&4,4,4\edge{dl} \\
&1,1&14\edge{urrrr}&12\edge{urrr}&&2,2&23\edge{u}&&3,3&34\edge{ulll}&&4,4& \\
&&&1\edge{ull}\edge{ul}\edge{u} &&2\edge{ull}\edge{u}\edge{ur} &&3\edge{ul}\edge{ur}\edge{urr} &&4\edge{urr}\edge{u}\edge{ulllllll}&&&\\
} } \caption {Constructing $L'$: L is a regular CW-complex
consisting of a $2$-cell, the square $1234$. $F(l)$ consists of
all multichains of rank $\leq r(l)$ in $(\hat{0},l]$ for all $l\in
L$. $L'\setminus\hat{0}$ is shown.}
\end{figure}
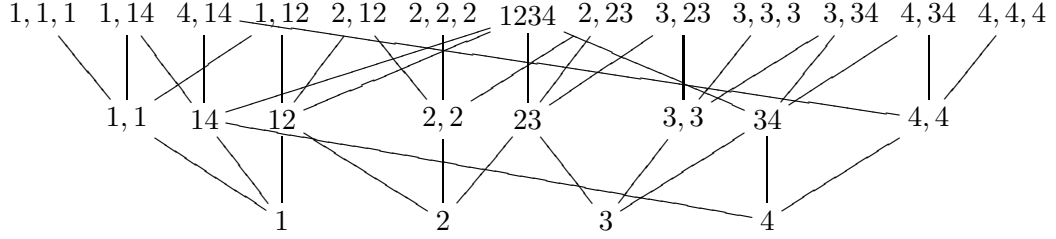

We define a new ranked meet semi-lattice $(L',<',r')$ as follows:
$L'=\cup_{l\in L\setminus \{\hat{0}\}}F(l)$ where the empty
multichain is the minimum $\hat{0}_{L'}$. In addition,
$r'(\underline{a})=\sum_{0\leq i\leq m}r(a_i)$ for
$\underline{a}\in L'$ as above, where the empty multichain has
rank $0$. We denote it in short by $L'$. See Figure \ref{FigL'}
for an illustration.

It is straightforward to verify that $L'$ is indeed a ranked meet
semi-lattice; we merely remark that for
$\underline{a},\underline{b}\in L'$ as above $\underline{a}\wedge
\underline{b}=(a_{\min (m,k)}\wedge b_{\min (m,k)}\leq...\leq
a_0\wedge b_0)$, which is an element of $L'$ as for $l\in L$ such
that $\underline{a}\in F(l)$ indeed $\underline{a}\wedge
\underline{b}\in F(l)$.
\begin{lem}\label{L' lemma}
Let $L$ be a ranked meet semi-lattice. If $L$ has the diamond
property then $L'$ has the parallelogram property.
\end{lem}

\begin{cor}\label{L' cor}
Let $L$ be a ranked meet semi-lattice, finite at every rank. If
$L$ has the diamond property then $L'$ satisfies Macaulay
inequalities (\ref{partial^{k-1}(n)}).
\end{cor}
$Proof$: This is immediate from Lemma \ref{L' lemma} and Theorem
\ref{THM}. $\square$

Before proving Lemma \ref{L' lemma}, let us mention that the $L'$
arising in this way include all posets of multicomplexes and all
meet semi-lattices with the diamond property. For the later, if
$L$ satisfies the diamond property, define
$F(l)=\{(l'):\hat{0}<l'\leq l\}\cup \{\emptyset\}$ for all $l\in
L$ to obtain $L'\cong L$. For a monomial $m$ in a multicomplex $M$
on the variables $x_1,...,x_n$, define $f(m)$ to be the unique
multichain of simplices obtained by dividing at each step by the
largest possible square free monomial, e.g. $m=x_1^5x_2x_4^3
\mapsto f(m)=(\{1\}\leq \{1\}\leq \{1,4\}\leq \{1,4\}\leq
\{1,2,4\})$. Denote by $\sigma(m)$ the largest simplex in the
multichain $f(m)$. $\sigma(m)=\{1,2,4\}$ in the example above. Let
$L=(\{\sigma(m): m\in M\}, \subset)$. It is (the face poset of) a
simplicial complex. For $\sigma\in L$ let $F(\sigma)=<f(m): m\in
M, \sigma(m)=\sigma>$ where $<>$ denotes the closure in the set of
multichains w.r.t. $\subset'$. Then $L'\cong M$ as ranked posets.
\newline \textbf{Remarks}:
(1) If $L$ is a regular CW-complex, $L'$ already gives us new
examples for which the inequalities (\ref{partial^{k-1}(n)}) hold,
see Figure 2.

(2) The construction $L\mapsto L'$ is a generalization of the
barycentric subdivision. If $L$ is a simplicial complex and $F(l)$
is the set of all chains (i.e. multichains $without$ repetitions)
in $(\hat{0},l]$ then $L'$ is the barycentric subdivision of $L$.
\\

$Proof\  of\  Lemma\  \ref{L' lemma}$: For every
$\underline{l'}\in L'$ consider the induced poset
$L'(\underline{l'})= \{\underline{y}\in L':\underline{l'}\leq'
\underline{y}\}$. An interval $[\underline{l'},\underline{l}]$,
$\underline{l}\neq \underline{l'}$, which is a chain in
$L'(\underline{l'})$ is of one of the following (intersecting) two
types: ($\prec '$ stands for the cover relation in $L'$.)

(1) $(\underline{l'}=\underline{l_0}\prec ' \underline{l_1}\prec
'...\prec ' \underline{l_m}=\underline{l})$ where there exists an
atom $u\in L$ such that for every $i$ where $1\leq i\leq m$
$\underline{l_i}$ is obtained from $\underline{l_{i-1}}$ by adding
$u$ to its lower end, denoted by
$\underline{l_i}=(u,\underline{l_{i-1}})$. In other words,
$[\underline{l'},\underline{l}]=\{\underline{l'}\prec '
(u,\underline{l'})\prec '(u,u,\underline{l'})\prec '...\prec '
(u,u,...,u,\underline{l'}) \}$.

(2) $(\underline{l'}\prec '\underline{l})$.

It follows from the fact that $L$ satisfies the diamond property
that indeed every interval not of type (1) nor of type (2) is not
a chain: let $\underline{a}\prec '\underline{b}\prec
'\underline{c}$ be a chain in such an interval, and assume by contradiction that it equals the interval $[\underline{a},\underline{c}]$. Combining this with the definition
of $L'$, we conclude that the multichains
$\underline{a},\underline{b},\underline{c}$ must have the same
length, i.e. same last index $m$ in the notation
$\underline{a}=(a_m\leq...\leq a_1\leq a_0)$. If $\underline{a}$
and $\underline{c}$ differ in at least two different indices,
denoted by $i$ and $j$, then clearly there are at least two
elements in the open interval $(\underline{a},\underline{c})$ -
just replace in $\underline{a}$ either $a_i$ with $c_i$ or $a_j$
with $c_j$. We are left to deal with the case where
$\underline{a}$ and $\underline{c}$ differ only in a single index,
$i$. As $\underline{a}\prec '\underline{b}\prec '\underline{c}$,
we conclude that $a_i\prec b_i\prec c_i$. By the diamond property
of $L$, there exists $d\in L$ such that $d\neq b_i$ and $a_i\prec
d\prec c_i$. Replacing $b_i$ with $d$ in $\underline{b}$ results
in a multichain $\underline{d}$ such that $\underline{a}\prec
'\underline{d}\prec '\underline{c}$; a contradiction. Thus indeed an interval not
of type (1) nor of type (2) is not a chain.

We now verify that $L'$ satisfies the parallelogram property.

Let $[\underline{l'},\underline{l}]$ be of type (1), and let
$\underline{x}\in (\underline{l'},\underline{l}]$,
$\underline{x}<'\underline{y}$, $(u,\underline{x})\nleq'
\underline{y}$. Then $\underline{x}=(u,\underline{x'})$ for some
multichain $\underline{x'}$. Let $d$ be the element in the
multichain $\underline{y}=(y_m\leq..\leq y_1\leq y_0)$ with the
same index as the index of $u$ at the lower end of $\underline{x}$
and let $c$ be the next indexed element in $\underline{y}$; put
$c=\hat{0}$ if $\underline{y}$ has the same last index as
$\underline{x}$. Then $u\not\leq c$. We will show now that there
exists $d'\in L$ such that $d'\prec d$, $c\leq d'$ and $u\nleq
d'$. Replacing $d$ with $d'$ in $\underline{y}$ we obtain a
multichain $\underline{y'}\in L'$ such that $\underline{y'}\prec '
\underline{y}$, $\underline{x'}<' \underline{y'}$ but
$\underline{x}\nleq' \underline{y'}$, as desired.

Let $\underline{\gamma}=(c<...<d)$ be a maximal chain in $[c,d]$
such that its element of minimal rank in its intersection with the
induced poset $L(c\vee u)$, denoted by $z$, is of maximal possible
rank. We need to show that $z=d$ (taking $d'$ as the element
covered by $d$ in $\underline{\gamma}$, we are done). Assume
$z\neq d$. Clearly $z\neq c$ (as $u\not\leq c$). Let $t\in
\underline{\gamma}$, $t\prec z$. By condition (\textbf{*}) of
Lemma \ref{cond* lemma}, there exists $t'\in L$ such that $t\prec
t'<d$ and $t'\neq z$. By the maximality of $z$, $t'\in L(c\vee
u)$. As $L$ is a meet semi-lattice, $c\vee u \leq t$,
contradicting the definition of $z$.

Let $[\underline{l'},\underline{l}]$ be of type (2), and not of
type (1). Let $\underline{l}<'\underline{y}$. By induction on the
rank $r(\underline{y})$ we will show the existence of
$\underline{y'}\in L'$ such that $\underline{y'}\prec '
\underline{y}$, $\underline{l'}<'\underline{y'}$ and
$\underline{l}\not\leq' \underline{y'}$. For
$r(\underline{y})=r(\underline{l})+1$, nonexistence of such
$\underline{y'}$ means that the chain $\underline{l'}\prec '
\underline{l}\prec ' \underline{y}$ is an interval, thus
$\underline{l}=(u,\underline{l'})$ for some atom $u\in L$, hence
$[\underline{l'},\underline{l}]$ is of type (1), a contradiction.
Thus, the case $r(\underline{y})=2$ is verified. Let
$\underline{t}\in [\underline{l},\underline{y}]$,
$\underline{t}\prec ' \underline{y}$. By induction hypothesis
there exists $\underline{z}$ such that
$\underline{l'}<'\underline{z}\prec ' \underline{t}$ and
$\underline{l}\not\leq' \underline{z}$. If the chain
$\underline{z}\prec '\underline{t}\prec '\underline{y}$ in $L'$ is
not an interval, let $\underline{y'}\in
(\underline{z},\underline{y})$, $\underline{y'}\neq
\underline{t}$. As $L'$ is a meet semi-lattice
$\underline{l}\not<' \underline{y'}$. We are left to deal with the
case $\underline{t}=(u,\underline{z})$ and
$\underline{y}=(u,\underline{t})$ for some atom $u\in L$. As
$\underline{l}\not\leq' \underline{z}$, the multichains
$\underline{l},\underline{t}$ have equal length, hence
$\underline{l}=(u,\underline{\tilde{l}})$ for some multichain
$\underline{\tilde{l}}$. As $[\underline{l'},\underline{l}]$ is
not of type (1), also $\underline{l'}=(u,\underline{\tilde{l'}})$
for some multichain $\underline{\tilde{l'}}$. Let us denote by
$\underline{\tilde{w}}$ the multichain obtained from
$\underline{w}$ by deleting its lower end $u$, where
$\underline{w}\in\{\underline{y},\underline{t},\underline{l},\underline{l'}\}$.
Looking at $L'(\underline{\tilde{l'}})$, by induction hypothesis
there exists $\underline{\tilde{y'}}\in
L'(\underline{\tilde{l'}})$ such that
$\underline{\tilde{l'}}<'\underline{\tilde{y'}}\prec '
\underline{\tilde{y}}$ and $\underline{\tilde{l}}\nleq'
\underline{\tilde{y'}}$. Then
$\underline{y'}=(u,\underline{\tilde{y'}})$ is as desired.
$\square$
\newline \textbf{Example T}: Let $T$ be a rooted tree such that
 all its leaves have the same distance $r$ from the root.
 Let $P(T)$ be the graded poset with $T$ as its Hesse diagram where
 the root is its maximal element.
 Add a minimum to $P(T)$ to obtain the ranked lattice $L(T)$.
 The parallelogram property trivially holds for $L(T)$,
 hence by Theorem \ref{THM} $L(T)$ satisfies Macaulay inequalities.
 (In this case, of course $f_0\geq f_1\geq...\geq f_{r}$,
 yet this family was not "trapped" by the previously known generalizations of
 Theorems \ref{KKthm} and \ref{Mthm}.)

\section{Face rings and algebraic shifting}\label{AlgShift}
\subsection{Shifting geometric meet semi-lattices}\label{GeomLattice}
We will associate an analogue of the exterior face ring to
geometric ranked meet semi-lattices, which coincides with the
usual construction for the case of simplicial complexes. Applying
an algebraic shifting operation, \`{a} la Kalai \cite{Kalai}, we
construct a canonically defined shifted simplicial complex, having
the same $f$-vector as its geometric meet semi-lattice.

Let $(L,<,r)$ be a ranked atomic meet semi-lattice with $L$ the set of its
elements, $<$ the partial order relation and $r:L\rightarrow
\mathbb{N}$ its rank function. We denote it in short by $L$. $L$ is called $geometric$ if
\begin{equation}\label{rank-ineq}
r(x\wedge y)+r(x\vee y)\leq r(x)+r(y)
\end{equation}
for every $x,y\in L$ such that $x\vee y$ exists. For example, the
intersections of a finite collection of hyperplanes in a vector
space form a geometric meet semi-lattice w.r.t. the reverse
inclusion order and the codimension rank. Face posets of
simplicial complexes are important examples of geometric meet
semi-lattices, where (\ref{rank-ineq}) holds with equality.

Adding a maximum to
a ranked meet semi-lattice makes it a lattice, denoted by $\hat{L}$, but the maximum may not have a rank.
Denote by $\hat{0},\hat{1}$ the minimum and maximum of $\hat{L}$, respectively, and by $L_i$ the set of rank $i$ elements in $L$. $r(\hat{0})=0$.

We now define the algebra $\bigwedge L$ over a field $k$ with
characteristic $2$. Let $V$ be a vector space over $k$ with basis
$\{e_u: u\in L_1\}$. Let $I_L=I_1+I_2+I_3$ be the ideal in the
exterior algebra $\bigwedge V$ defined as follows. Choose a total
ordering of $L_1$, and denote by $e_S$ the wedge product
$e_{s_1}\mathbf{\wedge}...\mathbf{\wedge} e_{s_{|S|}}$ where
$S=\{s_1<...<s_{|S|}\}$. Define:
\begin{equation}\label{I1-eq}
I_1=(e_S: S\subseteq L_1, \vee S=\hat{1}\in \hat{L}),
\end{equation}
\begin{equation}\label{I2-eq}
I_2=(e_S: S\subseteq L_1, \vee S\in L, r(\vee S)\neq |S|),
\end{equation}
\begin{equation}\label{I3-eq}
I_3=(e_S-e_T: T,S\subseteq L_1, \vee T=\vee S\in L, r(\vee
S)=|S|=|T|, S\neq T).
\end{equation}

(As ${\rm char} (k)=2$, $e_S-e_T$ is independent of the ordering
of the elements in $S$ and in $T$.) Let $\bigwedge L = \bigwedge
V/I_L$. As $I_L$ is generated by homogeneous elements, $\bigwedge
L$ inherits a grading from $\bigwedge V$. Let $f(\bigwedge
L)=(f_{-1},f_0,..)$ be its graded dimensions vector, i.e.
$f_{i-1}$ is the dimension of the degree $i$ component of
$\bigwedge L$.
\newline \textbf{Remark}:
If $L$ is the poset of a simplicial complex, then $I_L=I_1$ and $\bigwedge L$
is the classic exterior face ring of $L$, as in \cite{Kalai}.

The following proposition will be used for showing that $\bigwedge
L$ and $L$ have the same $f$-vector. Its easy proof by induction
on the rank is omitted.
\begin{prop}\label{rank}
Let $L$ be a geometric ranked meet semi-lattice. Let $l\in L$ and
let $S$ be a minimal set of atoms such that $\vee S=l$, i.e. if
$T\subsetneq S$ then $\vee T<l$. Then $r(l)=|S|$. $\square$
\end{prop}
\textbf{Remark}: The converse of Proposition \ref{rank} is also
true: Let $L$ be a ranked atomic meet semi-lattice such that every
$l\in L$ and every minimal set of atoms $S$ such that $\vee S=l$
satisfy $r(l)=|S|$. Then $L$ is geometric.

\begin{prop}\label{basis}
$f(\bigwedge L)=f(L)$.
\end{prop}
$Proof$:
Denote by $\tilde{w}$ the
projection of $w\in \bigwedge V$ on $\bigwedge L$.
We will show that
picking $S(l)$ such that $S(l)\subseteq L_1, \vee S(l)=l, |S(l)|=r(l)$
for each $l\in L$ gives a basis over $k$ of $\bigwedge L$, $E=\{\tilde{e}_{S(l)}: l\in L\}$.

As $\{\tilde{e}_{S}: S\subseteq L_1\}$ is a basis of $\bigwedge
V$, it is clear from the definition of $I_L$ that $E$ spans
$\bigwedge L$. To show that $E$ is independent, we will prove
first that the generators of $I_{L}$ as an ideal, that are
specified in (\ref{I2-eq}), (\ref{I1-eq}) and (\ref{I3-eq}),
actually span it as a vector space over $k$.

As $x\vee \hat{1}=\hat{1}$ for all $x\in L$, the generators of
$I_1$ that are specified in (\ref{I1-eq}) span it as a $k$-vector
space. Next, we show that the generators of $I_2$ and $I_1$ that
are specified in (\ref{I2-eq}) and in (\ref{I1-eq}) respectively,
span $I_1+I_2$ as a $k$-vector space: if $e_S$ is such a generator
of $I_2$ and $U\subseteq L_1$ then either $e_U\mathbf{\wedge}
e_S\in I_1$ (if $U\cap S\neq \emptyset$ or if $\vee(U\cup
S)=\hat{1}$) or else, by Proposition \ref{rank}, $r(\vee(U\cup
S))<|U\cup S|$ and hence $e_U\mathbf{\wedge} e_S$ is also such a
generator of $I_2$.

Let $e_S-e_T$ be a generator of $I_3$ as specified in
(\ref{I3-eq}) and let $U\subseteq L_1$. If $U\cap T\neq \emptyset$
then $e_T\mathbf{\wedge} e_U=0$ and $e_S\mathbf{\wedge} e_U$ is
either zero (if $U\cap S\neq \emptyset$) or else a generator of
$I_1+I_2$, by Proposition \ref{rank}; and similarly when $U\cap
S\neq \emptyset$. If $U\cap T= \emptyset = U\cap S$ then $\vee
(S\cup U)=\vee (T\cup U)$ and $|S\cup U|=|T\cup U|$. Hence, if
$e_S\mathbf{\wedge} e_U-e_T\mathbf{\wedge} e_U$ is not the obvious
difference of two generators of $I_1$ or of $I_2$ as specified in
(\ref{I1-eq}) and (\ref{I2-eq}), then it is a generator of $I_3$
as specified in (\ref{I3-eq}). We conclude that these generators
of $I_{L}$ as an ideal span it as a vector space over $k$.

Assume that $\sum_{l\in L}a_l \tilde{e}_{S(l)}=0$, i.e.
$\sum_{l\in L}a_l e_{S(l)}\in I_L$ where $a_l\in k$ for all $l\in
L$. By the discussion above, $\sum_{l\in L}a_l e_{S(l)}$ is in the
span (over $k$) of the generators of $I_3$ that are specified in
(\ref{I3-eq}). But for every $l\in L$ and every such generator $g$
of $I_3$, if $g=\sum\{b_S e_S: \vee S\in L, r(\vee S)=|S|\}$
($b_S\in k$ for all $S$) then $\sum\{b_S: \vee S=l\}=0$. Hence
$a_l=0$ for every $l\in L$. Thus $E$ is a basis of $\bigwedge L$,
hence $f(\bigwedge L)=f(L)$. $\square$

Now let us shift. Note that Kalai's algebraic shifting \cite{Kalai}, which was defined for the exterior face ring, can be applied to any graded exterior algebra finitely generated by degree $1$ elements. It results in a simplicial complex with an $f$-vector that is equal to the vector of graded dimensions of the algebra. This shows that any such graded algebra satisfies Kruskal-Katona inequalities! We apply this construction to $\bigwedge L$:

Let $B=\{b_u: u\in L_1\}$ be a basis of $V$. Then $\{\tilde{b}_S:
S\subseteq L_1\}$ spans $\bigwedge L$. Choosing a basis from this
set in the greedy way w.r.t. the lexicographic order $<_L$ on
equal sized sets ($S<T$ iff $\min (S\triangle T)\in S$), defines a
collection of sets:
$$\Delta_B(L)=\{S: \tilde{b}_S\notin {\rm span}_k\{\tilde{b}_T: |T|=|S|, T<_L S\}\}.$$
$\Delta_B(L)$ is a simplicial complex, and by Proposition
\ref{basis} $f(\Delta_B(L))=f(L)$. For a generic $B$,
$\Delta_B(L)$ is shifted. ($B$ is $generic$ if the entries of the
transition matrix form the standard basis to $B$ are algebraically
independent over a subfield of $k$. Alternatively, we can extend
$k$ by $n^2$ intermediates and consider the exterior algebra over
this bigger field, letting the transition matrix consist of those
intermediates. A collection of finite subsets of  $\mathbb{N}$,
$A$, is $shifted$ if $S\in A$ and $T$ that is componentwise not
greater than $S$ as ordered sets of equal size implies $T\in A$.)
Moreover, the construction is canonical, i.e. is independent both
of the chosen ordering of $L_1$ and of the generically chosen
basis $B$. It is also independent of the characteristic $2$ field
that we picked. We denote $\Delta(L)=\Delta_B(L)$ for a generic
$B$. For proofs of the above statements we refer to Bj\"{o}rner
and Kalai \cite{BK} (they proved for the case where $L$ is a
simplicial complex, but the proofs remain valid for any graded
exterior algebra finitely generated by degree $1$ elements).

We summarize the above discussion in the following theorem:
\begin{thm}\label{GL}
Let $L$ be a geometric meet semi-lattice, and let $k$ be a field
of characteristic $2$. There exists a canonically defined shifted
simplicial complex $\Delta(L)$ associated with $L$, with
$f(\Delta(L))=f(L)$. $\square$
\end{thm}
\textbf{Remarks}: (1) The fact that $L$ satisfies Kruskal-Katona
inequalities follows also without using our algebraic
construction, from the fact that it satisfies the diamond property
and applying Theorem \ref{Wthm}. The diamond property easily seen
to hold for all ranked atomic meet semi-lattices.

(2) A different operation, which does depend on the ordering of
$L_1$ and results in a simplicial complex with the same
$f$-vector, was described by Bj\"{o}rner \cite{Matroid
applications}, Chapter 7, Problem 7.25: totally order $L_1$. For
each $x\in L$ choose the lexicographically least subset
$S_x\subseteq L_1$ such that $\vee S_x=x$
($S_{\hat{0}}=\emptyset$). Define $\Delta_<(L)=\{S_x: x\in L\}$.
Then $\Delta_<(L)$ is a simplicial complex with the same
$f$-vector as $L$. An advantage in our operation is that it is
canonical (and results in a shifted simplicial complex). To see
that these two operations are indeed different, let $L$ be the
face poset of a simplicial complex. Then for any total ordering of
$L_1$, $\Delta_<(L)=L$. But if the simplicial complex is not
shifted (e.g. a $4$-cycle), then $\Delta(L)\neq L$.

\subsection{Shifting generalized multicomplexes}\label{MultiShift}
We will associate an analogue of the symmetric (Stanley-Reisner)
face ring with a common generalization of multicomplexes and
geometric meet semi-lattices. Applying an algebraic shifting
operation, we construct a multicomplex having the same $f$-vector
as the original poset.

Let $\mathbb{P}$ be the following family of posets: to construct
$P\in \mathbb{P}$ start with a geometric meet semi-lattice $L$.
Associate with each $l\in L$ the (square free) monomial
$m(l)=\prod_{a<l,a\in L_1}x_a$, and equip it with rank
$r(m(l))=r(l)$. Denote this collection of monomials by $M_0$. Now
repeat the following procedure finitely or countably many times to
construct $(M_0\subseteq M_1\subseteq...)$: Choose $m\in M_i$ and
$a\in L$ such that $x_a|m$, $\frac{x_a}{x_b}m\in M_i$ for all
$b\in L_1$ such that $x_b|m$, and $x_am\notin M_i$. $M_{i+1}$ is
obtained from $M_i$ by adding $x_am$, setting its rank to be
$r(x_am)=r(m)+1$ and let it cover all the elements
$\frac{x_a}{x_b}m$ where $b\in L_1$ such that $x_b|m$. Define
$P=\cup M_i$.

Note that the posets in $\mathbb{P}$ are ranked (not necessarily
atomic) meet semi-lattices with the parallelogram property, and
that $\mathbb{P}$ includes all multicomplexes (start with $L$, a
simplicial complex) and geometric meet semi-lattices ($P=M_0$).

For $P\in \mathbb{P}$ define the following analogue of the
Stanley-Reisner ring: Assume for a moment that $P$ is finite. Fix
a field $k$, and denote $P_1=\{1,..,n\}$. Let $A=k[x_1,..,x_n]$ be
a polynomial ring. For $j$ such that $1\leq j\leq n$ let $r_j$ be
the minimal integer number such that $x_j^{r_j+1}$ does not divide
any of the monomials $p\in P$. Note that each $i\in P$ of rank $1$
belongs to a unique maximal interval which is a chain; whose top
element is $x_i^{r_i}$. By abuse of notation, we identify the
elements in such intervals with their corresponding monomials in
$A$.

We add a maximum $\hat{1}$ to $P$ to obtain $\hat{P}$ and define the following ideals in $A$:

$I_0=(\prod_{i=1}^n x_i^{a_i}: \exists j\ 1\leq j\leq n,\
a_j>r_j),$

$I_1=(\prod_{i=1}^n x_i^{a_i}: \forall j\ a_j\leq r_j, \ \vee_{i=1}^n x_i^{a_i}=\hat{1}\in \hat{P}),$

$I_2=(\prod_{i=1}^n x_i^{a_i}: \vee_{i=1}^n x_i^{a_i}\in P, r(\vee_{i=1}^n x_i^{a_i})\neq\sum_i a_i),$

$I_3=(\prod_{i=1}^n x_i^{a_i}-\prod_{i=1}^n x_i^{b_i}: \vee_{i=1}^n x_i^{a_i}=\vee_{i=1}^n x_i^{b_i}\in P, r(\vee_{i=1}^n x_i^{a_i})=\sum_i a_i=\sum_i b_i),$

$I_P=I_0+I_1+I_2+I_3.$

Define $k[P]:=A/I_P$. As $I_P$ is homogeneous, $k[P]$ inherits a
grading from $A$. Let $f(k[P])=(f_{-1},f_0,..)$ where $f_i=\dim_k
\{m\in k[P]: r(m)=i+1 \}$ ($f_{-1}=1$).

The proof of the following proposition is similar to the proof of Proposition \ref{basis}, and is omitted.
\begin{prop}\label{basis2}
$f(k[P])=f(P)$. $\square$
\end{prop}
Denote by $\tilde{w}$ the projection of $w\in A$ on $k[P]$. Let
$B=\{y_1,..,y_n\}$ be a basis of $A_1$. Then
$$\Delta_B(P):=\{\prod_{i=1}^n y_i^{a_i}:
\prod_{i=1}^n\tilde{y_i}^{a_i} \notin {\rm
span}_k\{\prod_{i=1}^n\tilde{y_i}^{b_i}: \sum_{i=1}^n
a_i=\sum_{i=1}^n b_i, \prod_{i=1}^n y_i^{b_i}<_L \prod_{i=1}^n
y_i^{a_i} \}\}$$ is an order ideal of monomials with an $f$-vector
$f(P)$. (The lexicographic order on monomials of equal degree is
defined by $\prod_{i=1}^n y_i^{b_i}<_L \prod_{i=1}^n y_i^{a_i}$
iff there exists $j$ such that for all $1\leq t<j\ a_t=b_t$ and
$b_j>a_j$.) To prove this, we reproduce the argument of Stanley
for proving Macaulay's theorem (\cite{St}, Theorem 2.1): as the
projections of the elements in $\Delta_B(P)$ form a $k$-basis of
$k[P]$, then by Proposition \ref{basis2} $f(\Delta_B(P))=f(P)$. If
$m\notin \Delta_B(P)$ then $m=\sum\{a_n n: {\rm deg} (n)={\rm deg}
(m), n<_L m\}$, hence for any monomial $m'$ $m'm=\sum\{a_n m'n:
{\rm deg} (n)={\rm deg} (m), n<_L m\}$. But ${\rm deg} (m'm)={\rm
deg} (m'n)$ and $m'n<_L m'm$ for these $n$'s, hence $m'm\notin
\Delta_B(P)$, thus $\Delta_B(P)$ is an order ideal of monomials.
\newline \textbf{Remark}: For $B$ a generic basis the construction is canonical in
 the same sense as defined for the exterior case.

Combining Proposition \ref{basis2} with Theorem \ref{Mthm} we obtain
\begin{cor}\label{MP}
Every $P\in \mathbb{P}$ satisfies Macaulay inequalities
(\ref{partial^{k-1}(n)}). $\square$
\end{cor}
If $P$ is infinite, let $P_{\leq r}:=\{p\in P: r(p)\leq r\}$ and construct $\Delta(P_{\leq r})$ for each $r$.
Then $\Delta(P_{\leq r})\subseteq \Delta(P_{\leq r+1})$ for every $r$, and $\Delta(P):=\cup_r\Delta(P_{\leq r})$
is an order ideal of monomials with $f$-vector $f(P)$. Hence, Corollary \ref{MP} holds in this case too.

To conclude, I wish to address the following open question to the
readers:
\begin{prob}\label{Problem}
Find algebraic objects (such as standard graded rings) and notions of algebraic shifting that support Kruskal-Katona's and Macaulay's inequalities for the general combinatorial objects covered by Theorems \ref{Wthm} and \ref{THM}, respectively.
\end{prob}
\section*{Acknowledgments}
I deeply thank my advisor Prof. Gil Kalai for many helpful
discussions, and Prof. Anders Bj\"{o}rner for his comments on
earlier versions of this paper. Part of this work was done during
the author's stay at Institut Mittag-Leffler, supported by the ACE
network.


\begin{thebibliography}{99}
\bibitem{Matroid applications} A. Bj\"{o}rner, The homology and
shellability of matroids and geometric lattices, \textit{Matroid
applications} (N. White ed.), Cambridge Univ. Press, Cambridge,
1992.

\bibitem{BK} A. Bj\"{o}rner and G. Kalai, An extended Euler-Poincar\'{e} formula, \textit{Acta Math.}, \textbf{161} (1988), 279-303.

\bibitem{BV} A. Bj\"{o}rner and S. Vre\'{c}ica, On $f$-vectors and Betti numbers of multicomplexes, \textit{Combinatorica}, \textbf{17} (1997), 53-65.

\bibitem{Bol} B. Bollob\'{a}s, \textit{Combinatorics},
Cambridge Univ. Press, Cambridge 1986.

\bibitem{CL} G.F. Clements and B. Lindstr\"{o}m, A generalization of a combinatorial theorem of Macaulay, \textit{J. Combi. Th.},
\textbf{7}, (1969), 230-238.

\bibitem{EW} J. Eckhoff and G. Wegner, \"{U}ber einen Satz von Kruskal, \textit{Period. Math. Hung.},
\textbf{6}, (1975), 137-142.


\bibitem{Kalai} G. Kalai, A characterization of $f$-vectors of families of convex sets in $\mathbb{R}^{d}$, Part 1: Necessity of Eckhoff's conditions,
 \textit{Israel J. Math.}, \textbf{48}, (1984), 175-195.

\bibitem{M} F.S. Macaulay, Some properties of enumeration in the theory of modular systems, \textit{Proc. London Math. Soc.},
\textbf{26}, (1927), 531-555.



\bibitem{St} R.P. Stanley, Hilbert functions on graded algebras, \textit{Advances in Math.},
\textbf{28}, (1978), 57-83.

\bibitem{WW} D.L. Wang and P. Wang, Extremal configurations on a discrete torus and a generalization of the generalized Macaulay theorem, \textit{SIAM J. Applied Math.}, \textbf{33}, (1977), 55-59.

\bibitem{W} G. Wegner, Kruskal-Katona's theorem in generalized complexes,
\textit{Finite and Infinite Sets, Vol 2}, Coll. Math. Soc. \textbf{37}, (1984), 821-828.



\end{thebibliography}
\end{document}